\documentclass[letterpaper, 10 pt, conference]{ieeeconf}  

\IEEEoverridecommandlockouts                              
\usepackage{graphics} 
\usepackage{epsfig} 
\usepackage{amsmath} 
\usepackage{amssymb}
\usepackage{enumerate}
\usepackage{color}
\usepackage{float}
\usepackage{wrapfig}
\usepackage{graphicx}
\usepackage{caption}
\usepackage{subcaption}

\graphicspath{{figures/}}

\renewcommand{\th}{^{\rm th}}
\newcommand{\neighbI}{\mathcal{N}_{i}}
\newcommand{\neighb}{\mathcal{N}}

\newcommand{\vzero}{\mathcal{V}_0}

\newcommand{\al}{\alpha}

\newcommand{\R}{\mathbb{R}}
\newcommand{\Rn}{\mathbb{\R}^n}
\newcommand{\F}{\mathcal{F}}
\newcommand{\bg}{\mathcal{BG}}
\newtheorem{theo}{\bf Theorem}[section]

\newtheorem{corollary}[theo]{Corollary}

\newtheorem{definition}[theo]{Definition}
\newtheorem{example}[theo]{Example}

\newtheorem{lemma}[theo]{Lemma}

\newtheorem{problem}[theo]{\bfseries Problem}
\newtheorem{proposition}[theo]{Proposition}
\newtheorem{remark}[theo]{Remark}

\title{\LARGE \bf Formation Control with Pole Placement for Multi-agent Systems
}

\author{Ameer K. Mulla, Rachel K. Kalaimani, Debraj Chakraborty and Madhu N. Belur%
\thanks{A. K. Mulla, R. K. Kalaimani, D. Chakraborty and M. N. Belur are in the Department
of Electrical Engineering, Indian Institute of Technology Bombay,
India. \texttt{\small ameerkm@ee.iitb.ac.in, rachel@ee.iitb.ac.in, dc@ee.iitb.ac.in, belur@ee.iitb.ac.in}{\small {} }%
}}
\begin{document}

\maketitle
\thispagestyle{empty}
\pagestyle{empty}
\begin{abstract}
The problem of distributed controller synthesis for formation
control of multi-agent systems is considered. The agents
(single integrators) communicate over a communication graph and a
decentralized linear feedback structure is assumed. One of the agents is
designated as the leader. If the communication graph contains a 
directed spanning tree with the leader node as the root, then it is possible
to place the poles of the ensemble system with purely local feedback
controller gains. Given a desired formation, first one of the poles is
placed at the origin. Then it is shown that the inter-agent weights can be
independently adjusted to assign an eigenvector corresponding to the
formation positions, to the zero eigenvalue. Then, only the leader input is
enough to bring the agents to the desired formation and keep it there with
no further inputs. Moreover, given a formation, the computation of the
inter-agent weights that encode the formation information, can be calculated
in a decentralized fashion using only local information. 
\end{abstract}
\begin{keywords}
multi-agent systems, decentralized control, pole placement, formation control.
\end{keywords}
\section{Introduction}
Distributed control of multi-agent systems have been a popular research
topic with applications including sensor networks, rendezvous problems,
synchronization, surveillance, satellite formation (see e.g. \cite{Ren2004}
and the references therein) etc. In particular, consensus and formation
control for multi-agent systems communicating over a network, have
been an area of intense research since the initial formulation of
these problems in \cite{Fax2004,Olfati2004,Olfati-Saber2007}. In this framework, the individual
agents are required to reach a consensus or formation, i.e. fixed
relative position with respect to each other, with distributed control
and computation. Decentralized controllers which exploit the properties
of the graph Laplacian \cite{Fax2004, Jadbabaie2003, Olfati2004, Ren2005p}, to asymptotically bring all the
states to a specified consensus/formation, are particularly
well-studied. However, in most of the consensus/formation control
literature, the controller structure is pre-specified (to obtain the
closed loop in terms of the Laplacian) and hence the designer has
only limited freedom in specifying the rate of convergence of the
agents to the formation. In this article, we propose a distributed
controller synthesis method to achieve arbitrary pole placement for
formation control of multi-agent systems. 

On the other hand, moving formations are often required in many applications.
In such situations, the absolute location of the formation might change,
but the relative positions of the agents should remain the same. For
such a situation, \cite{Fax2001} suggested that the absolute desired
positions be communicated to each agent globally through a vector
of (possibly time varying) offsets. For moving the
agents in formation to a different location, these offsets need to
be globally recomputed and transmitted to each agent individually,
thereby making some part of the control essentially centralized. Existence
of a completely decentralized controller for moving formations was
characterized in \cite{Lafferriere2005}, but decentralized synthesis methods
for this problem seem to be unavailable. We address, in this paper,
formation synthesis with completely decentralized computation. We
show, that if each agent knows its desired relative position with
respect to its neighbors, then the formation can be moved to a new
absolute position, entirely by local re-computation of inter-agent
gains.\textcolor{red}{{} }

Assume there are $n$ agents indexed from $1$ to $n$ with identical
first order dynamics: 
\[
\dot{x}_i=u_i,\mbox{ for }i=1\mbox{ to }n.
\]
Each agent interacts with some of the other agents and their
communication pattern is depicted by a directed graph $G$, 
where the nodes of $G$ represent the agents. For a
given agent $i$, the set $\neighbI$ consists of the agents whose information
is available to $i$. We assume that the input to each agent is given
by a (completely general and decentralized) linear feedback law: 
\begin{eqnarray}
u_i & =\sum_{j\in\neighbI}\beta_{ij}x_j-\alpha_i x_i,
         \mbox{ for }i=1\mbox{ to }n-1.\label{ctrl:law-1}\\
u_n & =\sum_{j\in\neighb_n}\beta_{nj}x_j-\alpha_n x_n + v_n
\end{eqnarray}
Here $\alpha_i$ is the gain for the self feedback of each agent
and $\beta_{ij}$ is the gain for the feedback information from agent
$j$ to $i$. We designate the $n^th$ agent as the leader of the group and the
others as followers. It is assumed that an external input is allowed
on the leader, but not on any other agent. It is 
common \cite{Fax2004,Jadbabaie2003,Olfati2004,Ren2005p}
to assume, apriori special structure in $\beta_{ij}$ and $\alpha_i$,
namely assume $u_i(t) = -\sum_{j\in\neighbI} a_{ij} \left[x_i (t) - x_j(t)\right],\qquad i=1,...,n$.
For $a_{ij}=1$ the resulting closed loop dynamics is $\dot{x}(t)=-L_n x(t)$
where $x=[x_1 \enspace x_2\enspace...\enspace x_n ]^{T}$ and
$L_n \in\R^{n\times n}$ is the Laplacian (weighted for $a_{ij}$ not 
necessarily $1$) of the communication
topology. Clearly, the rate of convergence to consensus/formation is 
determined by the eigenvalues of the Laplacian matrix. In contrast, the values 
of $\alpha$ and $\beta$ in \eqref{ctrl:law-1} can be chosen
so as to arbitrarily place the closed loop poles of the ensemble of
systems. However, the computation of $\alpha$ and $\beta$ turns
out to be computationally hard for general pole placement. We prove,
that if $G$ contains a directed spanning tree 
(with $n\th$ agent as the root), and $G$ does not have any directed cycle the 
closed loop poles can be placed anywhere on
the real line just by adjusting the local feedback gains $\alpha$.
This allows us to decentrally compute
the self-feedback gains $\alpha$'s. The proof of this fact utilizes
properties of determinant expansions corresponding to bipartite graphs
\cite{lovasz, asratian}.

For formation control each agent is assumed
to know its relative position with respect to other neighboring agents.
We show that given any formation (say $\F=[f_1 ,...,f_n ],$
where $f_i $ is the absolute position of the $i$-th vehicle in
the desired formation), the gains $\alpha$ can be adjusted locally so as to 
bring the formation to $\F$
arbitrarily fast. Moreover, the agents remain in this formation even
after the external input to the leader is withdrawn. This is achieved
by placing a zero eigenvalue corresponding to the leader ($\alpha_n =0$)
and computing the inter-agent gains $\beta_{ij}$'s (with local information
only) so that $\F$ is the eigenvector corresponding to the
$0$ eigenvalue. Then all the agents can be forced into the formation $\F$ by 
just giving an external input $v_n$ corresponding to $f_n$ to the leader. For 
moving the entire formation to a new location
$\F^{'}=[f_1^{'},...,f_n^{'}]$, (preserving the same
relative position among the agents) only the leader (external) input
needs to be updated to $v_{n}^{'}$ (corresponding to $f_n^{'}$). The 
re-computation of the $\beta$'s
can be done again in a distributed fashion with local information
that propagates through the network.

Early versions of the consensus problem were presented in 
\cite{Chatterjee1977,Borkar1982}
but most of the recent research in this area is based on the 
theoretical
framework introduced for single integrator agents in 
\cite{Fax2004,Jadbabaie2003,Olfati2004,Ren2005p,Olfati-Saber2007}.
The consensus problem of the agents with second order dynamics was
studied in \cite{Hong2008,Olfati-Saber2006,Ren2008,Ren2007,Yu2010}.
In parallel, the formation control problem was introduced and developed
in \cite{Ren2004,Fax2004,Lafferriere2005}. 

The remaining article is arranged as follows. In section \ref{sec:prob}
the problem is formally defined and the main results of the article
are stated. In Section \ref{sec:prelim}, the preliminaries required
for proving the main results are given. The proofs of the main results
are given in \ref{sec:proofs}. The examples are given in Section
\ref{sec:exmp}.

\section{Problem Formulation and Main Results}\label{sec:prob}
We consider a system with $n$ agents as described in the introduction. We are 
specified with a particular formation $\F$ for the agents, i.e. the position 
to  which each of the agents should asymptotically reach from any given 
initial condition. Agent $n$ is considered as the leader for which the input 
$v_n$ is external. Our objective is to drive the agents to the desired 
formation using just $v_n$. It is expected that all other agents follow the 
leader. Hence all the agents are connected to the leader perhaps through other 
intermediate agents. This imposes that the graph $G$ has to have a rooted directed 
spanning tree. A spanning tree is a graph in which there are no cycles and
there is an undirected path between every two vertices. Of course, in
our case, the graph is directed. 
\begin{definition}
A directed spanning tree rooted at vertex $v$ is a spanning tree in which there is a 
directed path from $v$ to every other vertex in the graph.
\end{definition}
In this paper, the root is the vertex corresponding to the leader. 
Equation \eqref{ctrl:law-1} describes the input for each agent 
other than the leader. The closed loop dynamics are given by

\begin{equation}
\dot{x}=Ax+Bu_n \label{eq:system}
\end{equation}
 where  $$A=\left[\begin{array}{cccc}
-\alpha_1  & \beta_{12} & \cdots & \beta_{1n}\\
\beta_{21} & -\alpha_2  & \cdots & \beta_{2n}\\
\vdots & \vdots & \ddots & \vdots\\
\beta_{n-1,1} & \beta_{n-1,2} & \cdots & -\alpha_{n-1}\\
\beta_{n1} & \beta_{n2} & \cdots & \alpha_n
\end{array}\right], B=\left[\begin{array}{c}
0\\
0\\
\vdots\\
0\\
1
\end{array}\right]$$
Depending on $G$ the some of the $\beta$'s could be zero. We consider a specific example below to illustrate our assumption.

\begin{example}\label{exmp:prob}
For the case of $4$ agents and $1$ leader with interconnection defined by the graph in Figure \ref{fig:Comm-Graph},
the equations defining the agents' dynamics are as follows.
\begin{equation}\label{closed:loop}
\dot{x}=\begin{bmatrix}
\alpha_1  & 0 & \beta_{13} & \beta_{14} & 0\\
0 & \alpha_2  & 0 & 0 & \beta_{25}\\
0 & \beta_{32} & \alpha_3  & \beta_{34} & 0\\
0 & 0 & 0 & \alpha_4  & \beta_{45}\\
0 & 0 & 0 & 0 & \alpha_ 5
\end{bmatrix}x+\begin{bmatrix}0\\0\\0\\0\\1\end{bmatrix}v_n\end{equation}
\begin{figure}
\begin{center}
\includegraphics[width=0.2\textwidth]{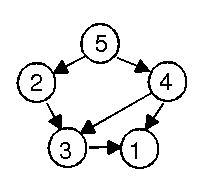}

\caption{Communication Graph\label{fig:Comm-Graph}}
\end{center}
\end{figure}
\end{example}
 
Next we choose some of the gains so as to make the eigenvalue corresponding to the leader, zero. We will show that if the leader pole is $0$ we benefit in the following aspects:
\begin{enumerate}
 \item After driving the leader to the desired position the external input is no longer required for achieving formation.
 \item It was discussed in the introduction that any given formation could be achieved using the inter-agent gains $\beta$'s. It turns out that the zero eigenvalue enables the decentralized computation of the $\beta$'s corresponding to any given formation.
\end{enumerate}

The zero eigenvalue can be achieved simply by letting $\alpha_n=0$ and $\beta_{nj}=0\enspace \forall j\in \neighb_n$.  This implies that the leader does not take any information from the other agents as well as it has no self-feedback. 
It is clear from Example \ref{exmp:prob} that the matrix $A$ has a pole at origin if $\al_5=0$.
The eigenvector, $\vzero$ corresponding to the $0$ eigenvalue determines the subspace to which the agents
converge. 
We show that it is possible to independently adjust the remaining $\beta$'s so that  such that $\vzero=\F$. Hence 
\begin{equation}\label{af=0}
 A\vzero=0
\end{equation}

The major question we address in this paper is about the ability to control
 the rate of convergence of each of the agents
to $\F$. In other words the question is whether we can place the poles of the closed loop system 
(except for the pole at origin) arbitrarily?
We assume that the set of (not necessarily distinct) closed loop poles 
are denoted as follows.
$$\Lambda:=\begin{bmatrix}\lambda_1,&\dots&,\lambda_{n-1},0\end{bmatrix}, \lambda_i\in\R\backslash 0.$$

Hence the feedback gains in addition to the constraints in \eqref{af=0} should also satisfy
\begin{equation}\label{cond:poles} \mbox { Roots of }\det(sI-A)= \Lambda.
 \end{equation}

 The following claims (to be made precise later) are the main contributions of this paper:
 \begin{enumerate}
  \item The poles can be placed arbitrarily by using only the self-feedback gains $\alpha$'s
  \item We use $u_n$ to drive the leader to $f_n$ and the rest of the agents follow appropriately to reach $\F$
  \item $\beta$'s can be computed in distributed manner so as to have $\vzero=\F$.
 \end{enumerate}
 
 In our problem the position of each agent is specified by a single dimension. 
In case each of the agent is assumed to move in a higher dimensional space, since the motion along each dimension is 
independent, it is enough to analyse the dynamics in one dimension.
A precise formulation of the problem is as follows.
\begin{problem}
Assume $n$ agents, with $n\th$ agent as leader, with dynamics $\dot{x}=u$ whose interconnection
is described by a directed connected graph $G$. Except the leader, the control law for each of the agents 
is given by \eqref{ctrl:law-1}. Let the final formation to which the agents should reach 
be specified by the vector $\F \in\Rn$.  Let $\Lambda$ denote the set of specified poles
 of the system.
\begin{enumerate}
 \item Find sufficient condition on $G$ which enables calculation of the feedback gains $\alpha_i$ and $\beta_{ij}$ 
in order that the agents converge to $\F$ with the convergence rate given by the set $\Lambda$. 
\item When are the gains $\alpha_i$ and $\beta_{ij}$ in the control law unique?
\item What is the minimum interaction, i.e. the minimum number of edges in $G$ required
 to achieve statement $1$?

\end{enumerate}
\end{problem}
\subsection{Main Results}
We state the main results of the paper here and the proofs are in the following sections.
The solution to our problem  involves finding the feedback gains satisfying constraints \eqref{af=0} and
\eqref{cond:poles}.  The following theorem provides a sufficient condition on $G$ which ensures precise
calculation of the feedback gains from the above mentioned constraints.
\begin{theo}\label{formation:rate}
 Assume there are $n$ agents with one leader with the dynamics of each agent as $\dot{x}_i=u_i$.
 Let a directed graph $G$ represent the interconnection
between them which defines the structure of the control laws as given in \eqref{ctrl:law-1}. 
The dynamics 
of the closed loop is $\dot{x}=Ax+e_n u_n$, with $u_n$ as the external input to the leader.
Let $\F\in\Rn$ denote the final formation of the
agents and suppose $\Lambda$ denotes the desired set of closed loop poles of the system. 
The following conditions are sufficient to ensure that the agents converge to $\F$ with 
the rates specified by
$\Lambda$.
\begin{enumerate}
 \item $G$ contains a directed spanning tree with vertex $n$ corresponding to leader as root.
\item $G$ does not have a directed cycle.
\end{enumerate}
A controller having the decentralized controller structure proposed in \eqref{ctrl:law-1} that achieves this is the solution
 to the following system of equations.
\begin{eqnarray}
&\alpha_i=\lambda_i.\label{alpha:unique}\\
&f_i\lambda_i+\sum_{j\in \neighbI}\beta_{ij}f_j=0, \mbox{ for } i=1:n-1. \label{beta:notunique}
\end{eqnarray}
\end{theo}

The absence of directed cycles in  $G$ allows just the
`inner' loops within each agent to decide the closed loop poles
of the interconnected system;
this is precisely the reason that arbitrary {\em real} poles can be
placed.
On the other hand, $\beta_{ij}$ plays a role in just obtaining the desired
formation vector. In this way, the steady-state formation configuration and
closed loop pole placement exhibit the {\em separation principle}.


In the above theorem the self feedback gains were unique (upto ordering of the closed loop poles),
 which were precisely the poles specified. 
The next question is about the uniqueness of $\beta$'s which is addressed in the following corollary.
It is shown that this case occurs when the interaction between the agents is minimum.


\begin{corollary}\label{min:inter}
The minimum number of directed interaction required between $n$ agents 
for converging to a specified formation $\F$ with the specified rates given by $\Lambda$ is $n-1$.
In this case the feedback gains $\alpha$ and $\beta$ are unique which are given by
$$\alpha_i=\lambda_i \mbox{ and }
\beta_{ij}=-\frac{f_i\alpha_i}{f_j},j\in \neighbI  \mbox{ for } i=1:n-1.$$
\end{corollary}

The third claim about the distributed computation of $\beta$'s is explained in the following subsection.
\subsection{Calculation of $\beta$'s}\label{subsec:rmk}
 When the formation is specified as the desired relative states of the agents, the leader position dictates the absolute state values of the agents in the formation $\F$. We assume that the absolute positions $f_1 ,...,f_n $ are not known to the agents $1,...,n-1$, but the relative position information i.e.  $\gamma_{ij}=f_i -f_j \enspace\forall j\in\neighbI $ is known to the agents. In addition to $\gamma_{ij}$'s the information about $f_j \enspace\forall j\in\neighbI $
is available through the communication channel.
The $i\th$ agent should converge to $f_i $ corresponding to the
leader's target position $f_n $ with given rate of convergence.
The absolute value of formation target for $i\th$ agent is specified
as $f_i =f_j +\gamma_{ij}\enspace\forall j\in\neighbI $.
Note that, $\gamma$'s should satisfy $\gamma_{ij}+\gamma_{jk}=\gamma_{ik}$ for any $i,j,k$ where $i,j,k$ refer to the indices of any three agents.

Every time the target position for the
leader $f_n $ is changed, the corresponding $\F$ changes but $\gamma_{ij}$'s are fixed.
Every $\F$ changes, we need to find the $\beta_{ij}$'s
so that $A\F=0$. From $A$ defined in equation (\ref{eq:system}),
we can write 
\[
-\alpha_i f_i +\sum_{j\in\neighbI }\beta_{ij}f_j =0
\]
 The $\alpha_i $ and $\gamma_{ij}$ are fixed and information about
$f_j $ is available through the communication channel.  $\beta_{ij}$ appear only in the equation corresponding
to the $i\th$ row of $A\F=0$. Since this equation involves
other parameters that are available to the $i\th$ agent locally,
the $\beta_{ij}$ can be calculated locally.

\section{Preliminaries}\label{sec:prelim}
A directed graph $G(V,E)$ refers to a set of vertices $V$ and a set of edges $E$ in which each element is
an ordered pair of vertices. A path in a directed graph $G$ is a finite sequence of vertices 
$v_0,v_1 ,\dots,v_n $ such that $(v_i,v_{i+1})\in E(G)$ for $i=0:n-1$. If $v_n=v_0$ then 
this path is a cycle in the directed graph $G$. We next consider bipartite graphs which are not directed.
\subsection{Bipartite graphs} \label{bi:graph}
A graph $\bg=(V,E)$ with vertex set $V$ and edge set $E$ is said 
to be bipartite  if $V$ can be partitioned into two subsets $V_1$ and $V_2$ such
that 
no two vertices from the same subset are adjacent. If we assign weights for the edges then we have a
weighted bipartite graph.
A set of edges $M$ in a graph $\bg$ is called a \emph{matching} if no two edges in
$M$ are adjacent. 
A maximal matching is a matching with maximum number of edges. A graph can have
more
than one maximal matching. The cardinality of a 
matching $M$, denoted by $|M|$, is defined as the number of edges in $M$. 
 $\bg$ is said to have a perfect matching if $|V_1|=|V_2|$ and there exists a matching $M$ such that 
$|M|=|V_1|=|V_2|$.

We show how a square polynomial matrix $P(s)\in\R^{n \times n}[s]$ can be associated to a  weighted bipartite graph.
The set $R$ and $C$ denote the rows
and columns 
of the polynomial matrix and are the two disjoint vertex sets of the bipartite
graph $\bg$, i.e. $|R|=|C|=n$.
By definition an edge exists in the bipartite graph between 
vertex $v_i \in R$ and $v_j \in C$  if the $(i,j)$-th entry of $P(s)$ is
nonzero. The weight of the edge is the corresponding entry itself.

Next we describe the relation between the 
determinant of $P$  and perfect matchings of the bipartite graph $\bg$ associated
to $P$.  Let $M$ be a perfect matching in $G$. Then the \emph{product} of the weights of all edges in $M$
 corresponds to a nonzero term in 
the determinant expansion of $P$. 
The determinant expansion of $P$ is the sum over all perfect matchings 
in $G$ (with suitable signs).
See \cite{bab:fra}. A cycle in a bipartite graph is defined similarly as in directed graphs except that the edges in this
case are not directed. We define an alternating cycle.
\begin{definition}
 A cycle in bipartite graph, $\bg$ is said to be \emph{alternating relative} to a matching $M$ if its edges 
are alternately in $E(\bg)\backslash M$ and $M$
\end{definition}
The following proposition is about a condition for the existence of more than one matching in $\bg$.
\begin{proposition}\cite{asratian}\label{altr:cycle}
 Let bipartite graph, $\bg$ have a perfect matching $M$. Then every other perfect matching can be obtained from
$M$ by a sequence of transfers along alternating cycles relative to $M$.
\end{proposition}
The above proposition also implies that if there exist more than one perfect matching
in $\bg$ then there exists a cycle in $\bg$.

\section{Proof of Main Results}\label{sec:proofs}
In this section we prove the main results of this paper. We begin with the proof of our first main result
 Theorem \ref{formation:rate}. Since poles are the roots of the determinant of $(sI-A)$ we use
 the relation between the determinant of a matrix and matchings in the
corresponding bipartite
graph. In this regard we state the following lemma which relates the existence of a cycle in $G$ and corresponding $\bg$.

\begin{lemma}\label{cycle_g}
 Let $G$, a directed graph represent the interconnection among $n$ agents. 
The feedback laws with respect to $G$ are as in \eqref{ctrl:law-1} and let the dynamics 
of the closed loop system be $\dot{x}=Ax+e_n u_n$. Let $\bg$ denote the bipartite graph constructed for $(sI-A)$.
Then the undirected bipartite graph $\bg$ has a cycle if and only if the directed graph $G$ has a directed cycle.
\end{lemma}
\begin{proof}
 Assume in $\bg$, the vertex set $R$ and $C$ correspond respectively to the rows and columns of 
$sI-A$ as explained in Section \ref{bi:graph}. Let $u_i$ and $v_i$ for $i=1:n$ represent the vertices in
$R$ and $C$ respectively. Hence an edge between $u_i$ and $v_i$ correspond to the diagonal terms in $(sI-A)$.
An edge between $u_i$ and $v_j$, $i\neq j$ implies that in $G$, $j\in \neighbI $. This means that agent $i$
takes data from agent $j$. We assume 
$\bg$ has a cycle and prove $G$ also has a cycle. Assume we index the vertices in
$R$ and $C$ suitably such that the cycle in $\bg$ is $u_1v_1u_2v_2,\dots,u_kv_ku_1$. 
Then $2\in \neighb_1,\dots,k\in \neighb_k-1, 1\in \neighb_k$. This implies that there is a cycle in $G$.
The proof that $G$ has a cycle implies $\bg$ has a cycle is straightforward and hence skipped.
\end{proof}
{\bf Proof of Theorem \ref{formation:rate}:} \\
Since $G$ has a directed spanning  tree with the leader at the root, all the agents can follow the leader.
 The condition for placing the poles at the given set $\Lambda$ is 
$$\mbox { Roots of }\det(sI-A)= \Lambda. $$
Let $\bg$ represent the bipartite graph constructed for $(sI-A)$. 
Then because of $sI$, $\bg$ has at least one perfect matching. $G$ does not have a cycle. Hence from Lemma
\ref{cycle_g} $\bg$ does not have a cycle. Hence from Proposition \ref{altr:cycle} we conclude that there 
exist only one perfect matching. Hence as mentioned in Section \ref{bi:graph}, 
the determinant of $sI-A$ is given by the product
$s\prod_{i=1}^{n-1}(s-\alpha_i)$. Therefore given the set of poles, $\Lambda$ we assign $\alpha_i=\lambda_i$,
for $i=1:n-1$.
In order that the agents reach the formation specified by
$\F$ asymptotically, it is required that the eigenvector corresponding to $0$ eigenvalue should be $\F$.
Hence we have the condition $A\F=0$ which is same as
\begin{equation}
f_i\alpha_i+\sum_{j\in \neighbI }\beta_{ij}f_j=0, \mbox{ for } i=1:n-1.\end{equation}
Substituting $\alpha_i=\lambda_i$ we get
\begin{equation}\label{eq1}
f_i\lambda_i+\sum_{j\in \neighbI }\beta_{ij}f_j=0, \mbox{ for } i=1:n-1.\end{equation}
Since $G$ has a  tree, $|\neighbI |\geqslant 1$ for $i=1:n-1$. Hence from \eqref{eq1},
the map from  $\beta_{ij}$'s to $\F$ is a linear surjective map. Hence we can calculate the
required gains.
$\hfill \blacksquare $\\
\begin{remark}
 Ability to place all the poles amounts to surjectivity
(see \cite{RosenthalWillems}) of
a certain map from the space of all parameters: $\alpha_i$ and
$\beta_{ij}$ to
the space of all coefficients of possible closed loop characteristic
polynomials.
This map turns out to be straightforward when $G$ has no directed cycle, but
this map is a {\em polynomial} map in both types of
indeterminates $\alpha_i$ and $\beta_{ij}$ in the presence of a cycle.
A linearization of this map is done in
\cite{RosenthalWillems}.
Proving surjectivity of the linearized map within the decentralized
controller structure and proving that a directed cycle can indeed result
in non real closed loop poles requires further investigation along these
lines.
\end{remark}

{\bf Proof of Corollary \ref{min:inter}:}\\
 Since $G$ has a  tree the number of edges in $G$ is at least $n-1$. 
We show that with $G$ having $n-1$ edges it is possible to make the agents 
converge to a specified formation with the desired rate. If $G$ has $n-1$ edges then 
$G$ is just a  tree and hence there are no cycles.
Therefore from Theorem \ref{formation:rate} we have
$\alpha_i=\lambda_i$ for $i=1:n-1$. Equation \eqref{beta:notunique} in the same theorem is considered.
$$f_i\lambda_i+\sum_{j\in \neighbI }\beta_{ij}f_j=0, j\in \neighbI , \mbox{ for } i=1:n-1.$$
Since $|\neighbI |=1$, we have unique $\beta$'s which is given by
$$\beta_{ij}=-\frac{f_i\alpha_i}{f_j},j\in \neighbI ,\mbox{ for } i=1:n-1.$$ 
Thus the minimum number of interaction required among $n$ agents for converging to $\F$ with the rates given 
by the set $\Lambda$ is $n-1$.
$\hfill\blacksquare$\\
The above corollary implies that the non uniqueness of $\beta$ arises when $G$ has 
more than $n-1$ edges. However $\alpha$ is still unique. 
Therefore the next possible 
question is about the maximum number of edges $G$ can allow without losing the uniqueness of $\alpha$
which is explained in the following remark.  This also
gives the maximum interaction that is possible while retaining the ability to  calculate the feedback gains.

\begin{remark}\label{max:inter}
 The maximum number of directed interaction allowed between $n$ agents 
  while converging to $\F$ with the specified rates given by $\Lambda$ is $n(n-1)/2$.
When $\bg$ does not have a cycle the vertices in $R$ and $C$ can be indexed suitably such that the 
corresponding matrix is upper triangular. In $(sI-A)$ the diagonal terms correspond to the self feedback gains. 
The rest of the entries indicate the presence of a communication link with other agents. Hence there
can at most be $n(n-1)/2$ entries in $(sI-A)$ in order to retain the upper triangular property. Therefore the
 maximum allowed communication links such that we can calculate the feedback gains is $n(n-1)/2$.
\end{remark}

\section{Examples}\label{sec:exmp}

\begin{example}\label{exmp:1}

Consider a system with five agents. Let agent $5$ be the leader. The communication graph for the system is given in Figure \ref{fig:Comm-Graph}. Notice that the communication graph has 6 edges i.e. a spanning tree and two additional edges but no cycles. 

With the decentralized control law (\ref{ctrl:law-1}) for the followers, 
the resultant system dynamics is given by,
\[
\dot{x}=Ax+Bu
\] with
\[
A=\begin{bmatrix}
\alpha_1 & 0 & \beta_{13} & \beta_{14} & 0\\
0 & \alpha_2 & 0 & 0 & \beta_{25}\\
0 & \beta_{32} & \alpha_3 & \beta_{34} & 0\\
0 & 0 & 0 & \alpha_4 & \beta_{45}\\
0 & 0 & 0 & 0 & \alpha_5\\     
\end{bmatrix}\mbox{ and }B=\begin{bmatrix}0\\0\\0\\0\\1\end{bmatrix}
\]
A pole at origin is placed by making $\alpha_5=0$. The remaining system poles are to be placed $-3, \enspace-3.5, \enspace-4 \mbox{ and } -5$. By Theorem \ref{formation:rate}, we know that the values of $\alpha_i$'s are indeed system poles. Thus we get unique (upto the ordering) values of $\alpha_i$'s as $\alpha_1=-3, \alpha_2=-3.5, \alpha_3=-4 \mbox{ and } \alpha_4=-5$. The $\beta$'s are chosen for achieving formation $\F=[-3 \enspace 2 \enspace -2 \enspace -1 \enspace 1]^T$. Figure \ref{fig:tree+form} demonstrates that for given $\F$, the $\beta$'s can be non unique.

\begin{figure}
        \centering
        \begin{subfigure}[b]{0.5\textwidth}
                \centering
                \includegraphics[width=\textwidth]{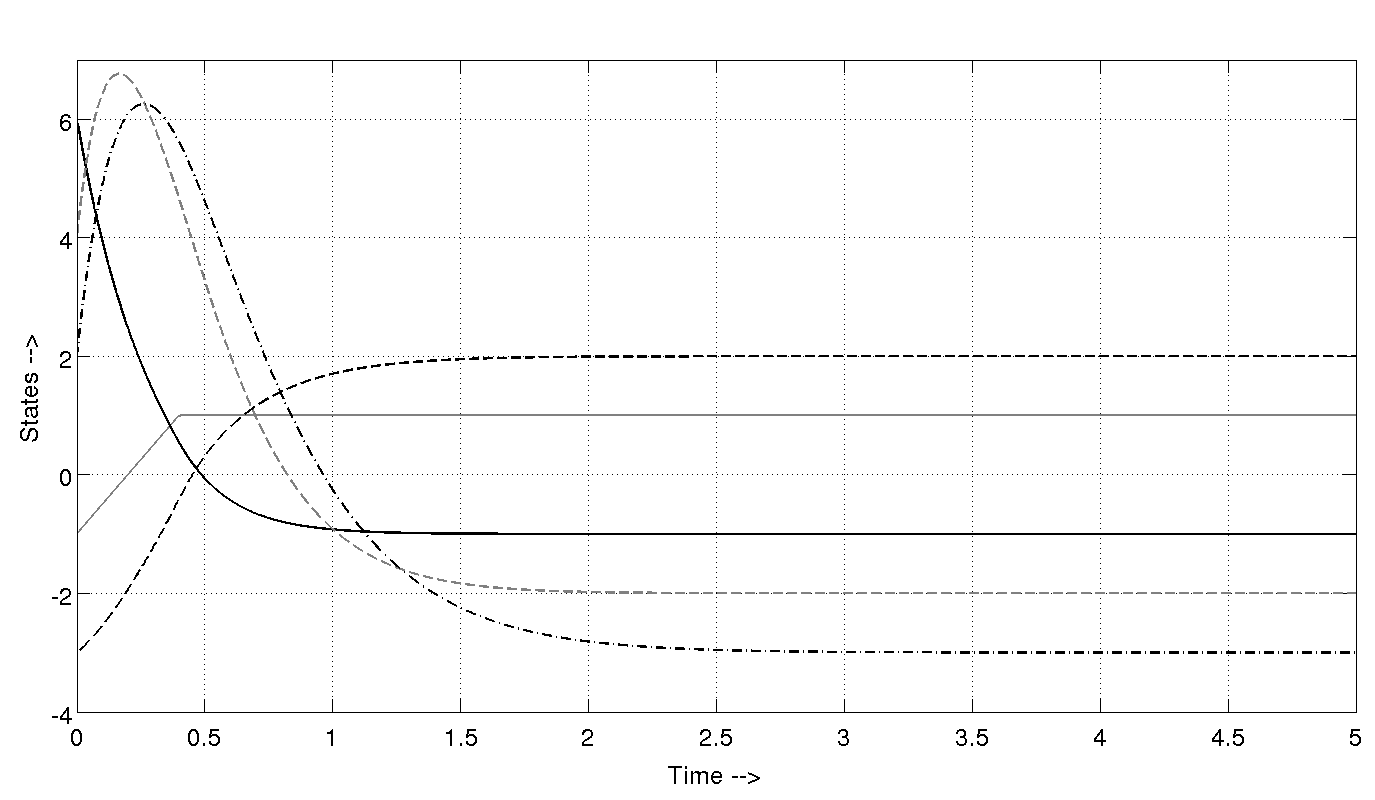}
                \caption{With $\beta_{13}=1$, $\beta_{14}=2$, $\beta_{32}=2$ and $\beta_{34}=2$}
                \label{subfig:tree+cons1}
        \end{subfigure}%
        ~ 

        \begin{subfigure}[b]{0.5\textwidth}
                \centering
                \includegraphics[width=\textwidth]{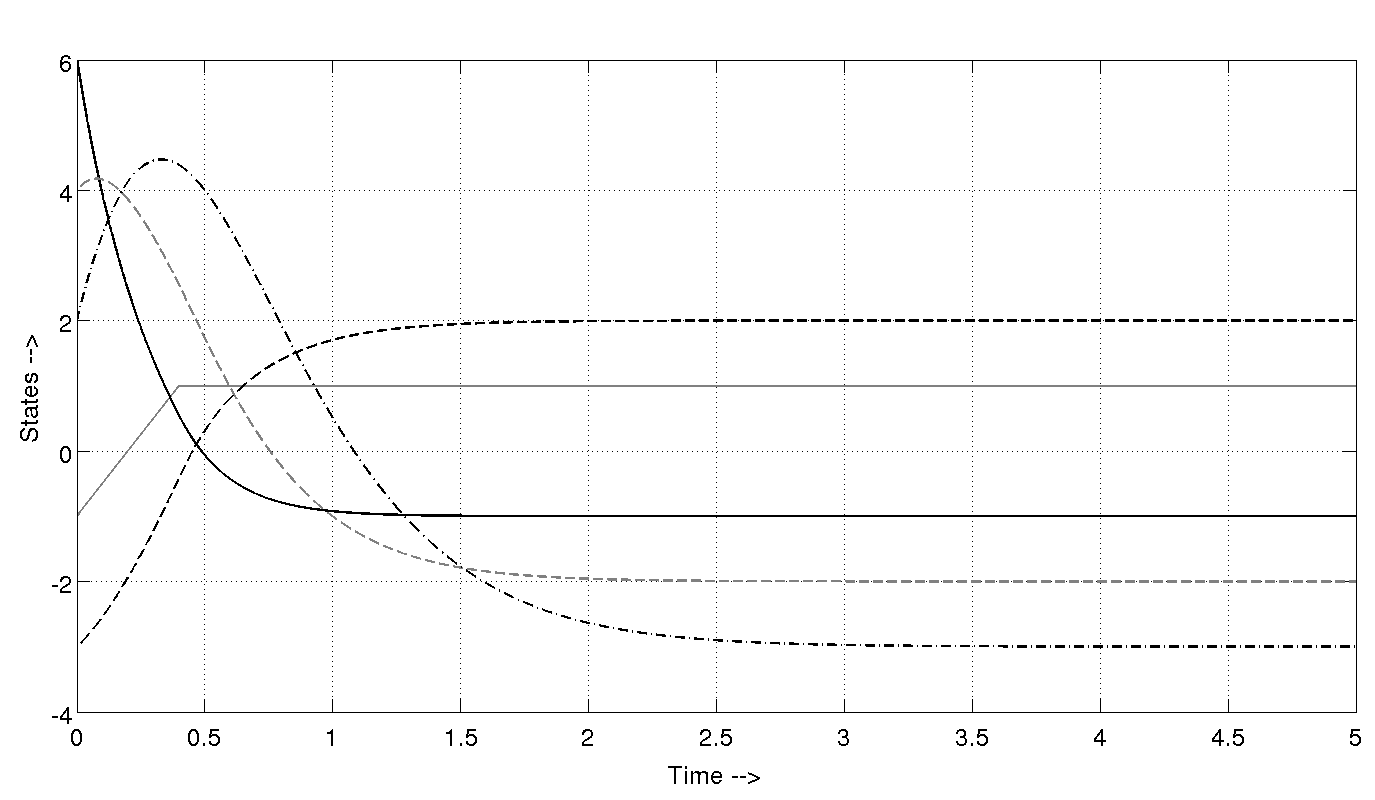}
                \caption{With $\beta_{13}=4$, $\beta_{14}=-1$, $\beta_{32}=1$ and $\beta_{34}=3$}
                \label{subfig:tree+cons2}
        \end{subfigure}
        ~ 
                \caption{Formation $\F$ with different $\beta$'s but same eigenvalues}\label{fig:tree+form}
\end{figure}

\end{example}

\begin{example}\label{exmp:2}
 This example demonstrates the non uniqueness of $\alpha_i$'s for pole placement and consensus when the communication topology has a cycle. Consider a simple system with three agents. Let agent $3$ be the leader. The communication graph for the system is given in Figure \ref{fig:graph-exmp2}. Notice that the communication graph has 3 edges i.e. a spanning tree and one additional edges out of which two edges form a directed cycle. 
\begin{figure}
  \begin{center}
    \includegraphics[width=0.2\textwidth]{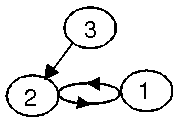}
  \end{center}
  \caption{Communication graph for example \ref{exmp:2}}\label{fig:graph-exmp2} 
\end{figure}

 With the decentralized control law (\ref{ctrl:law-1}) for the followers, we have,

\[
\dot{x}=\begin{bmatrix}
\alpha_1 & \beta_{12} & 0\\
\beta_{21} & \alpha_2 & \beta_{23}\\
0 & 0 & 0
\end{bmatrix}+
  \begin{bmatrix}0\\0\\1\end{bmatrix}  u_3
\]
$\alpha_3$ is set to $0$ so as to achieve a pole at origin. The system poles are no longer independent of $\beta$'s.
For achieving consensus i.e. $\F=[1\enspace1\enspace...\enspace1]^T$ with the agents converging to $f_i=4 \enspace \forall i$ and placing the remaining system poles at $-4 \mbox{ and } -5$, we can have different values of $\alpha$'s and $\beta$'s. Figure \ref{fig:cycle_cons} demonstrates that with different $\alpha$'s and $\beta$'s, the consensus can be reached.

\begin{figure}
        \centering
        \begin{subfigure}[b]{0.5\textwidth}
                \centering
                \includegraphics[width=\textwidth]{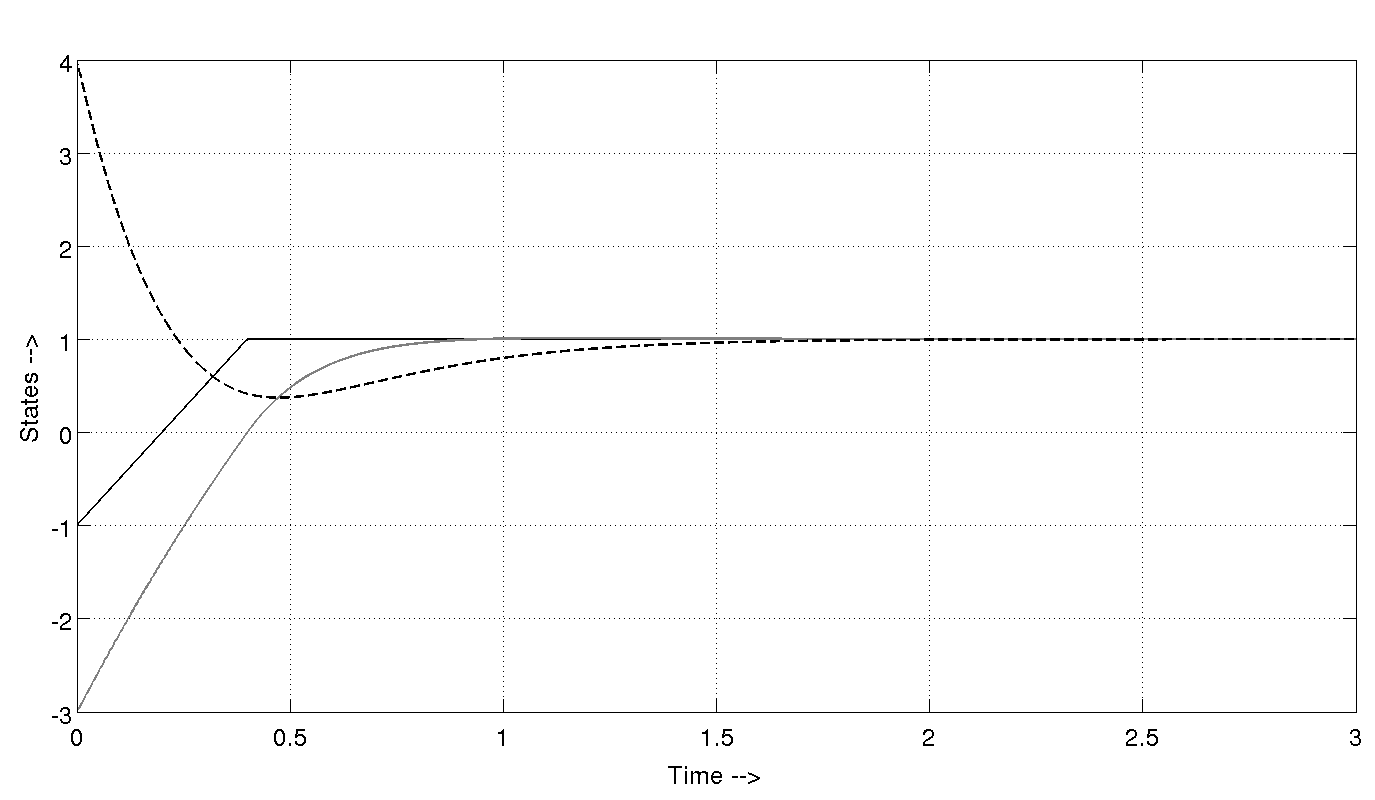}
                \caption{With $\alpha_1 =3$, $\alpha_2 =6$, $\beta_{12}=3$, $\beta_{21}=-2/3$ and $\beta_{23}=20/3$}
                \label{subfig:cycle_cons1}
        \end{subfigure}%
        ~ 

        \begin{subfigure}[b]{0.5\textwidth}
                \centering
                \includegraphics[width=\textwidth]{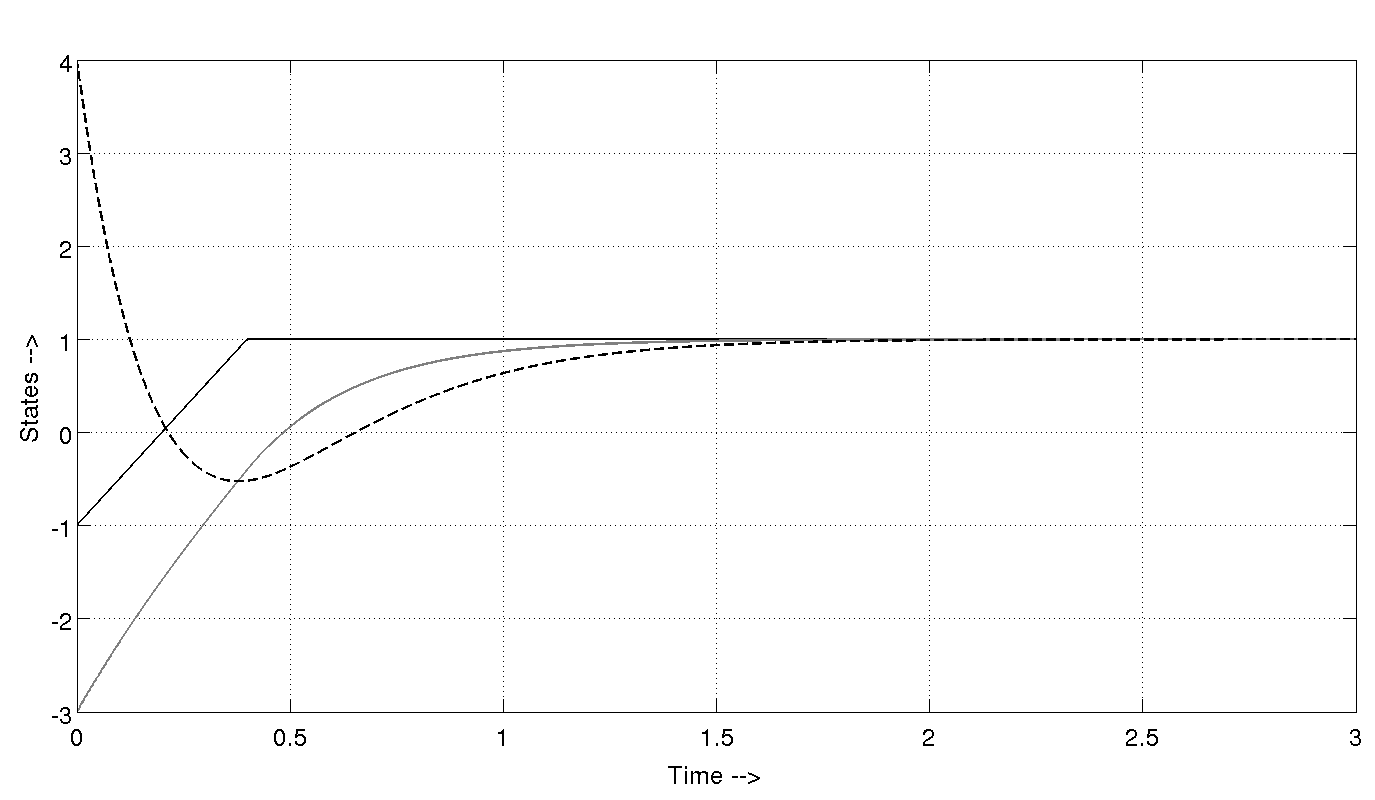}
                \caption{With $\alpha_1 =5$, $\alpha_2 =4$, $\beta_{12}=5$, $\beta_{21}=0$ and $\beta_{23}=4$}
                \label{subfig:cycle_cons2}
        \end{subfigure}
        ~ 
                \caption{Consensus $\F$ with different $\alpha$'s and $\beta$'s but same eigenvalues}\label{fig:cycle_cons}
\end{figure}
The eigenvalues also remain the same for these values of $\alpha_i$'s and $\beta_{ij}$'s.
\end{example}

If the communication graph $G$ has a directed cycle, the system poles are different from the self-feedback gains ($\alpha$'s) of the agents. So we can also have complex poles for the closed loop system. The following example demonstrates that complex poles can be placed with real $\alpha$'s and real $\beta$'s.
\begin{example}\label{exmp:3}
  
 Consider a communication graph in Figure \ref{fig:graph-exmp3}. 
 \begin{figure}
  \begin{center}
    \includegraphics[width=0.2\textwidth]{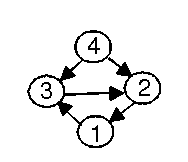}
  \end{center}
  \caption{Communication graph for example \ref{exmp:3}}\label{fig:graph-exmp3} 
\end{figure}
 The $\alpha$'s and $\beta$'s are chosen so as to achieve formation $\F=[2\enspace-1\enspace-2\enspace1]^T$ with system poles at $-5.5377$, $-2.7312 + 1.5140i$ and $-2.7312 - 1.5140i$. The closed loop system is given by 
\[
\dot{x}=\begin{bmatrix}
           -3 & 1 & 0 & 7\\
           0 & -4 & 2 & 0\\
           -3 & 0 & -4 & -2\\
           0 & 0 & 0 & 0                
          \end{bmatrix}x
  +\begin{bmatrix}0\\0\\0\\1\end{bmatrix}  u_n  .
\] 
The formation is shown in Figure \ref{fig:cycle_form}
\begin{figure}
  \begin{center}
    \includegraphics[width=0.5\textwidth]{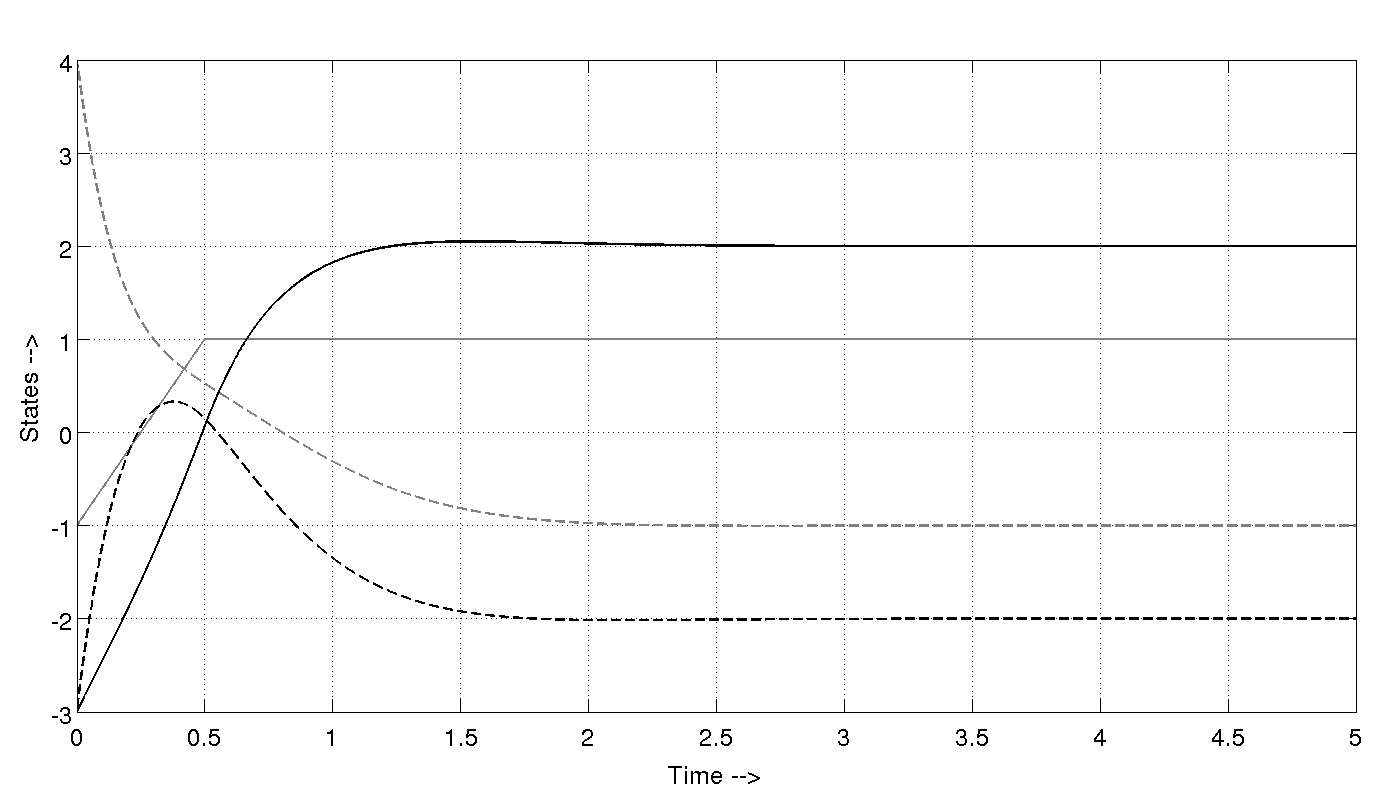}
  \end{center}
  \caption{Formation for the system in Example \ref{exmp:3}}\label{fig:cycle_form} 
\end{figure}    

\end{example}
\begin{example}\label{exmp:4}
 In this example we consider a network of agents in two dimensional state space. Motion in each dimension is governed by separate input. So the agents are still single integrators. We consider a system with $6$ agents communicating to each other through graph $G$ shown in figure \ref{fig:graph-tree}. Agent $6$ is the leader. So $\alpha_6=0$ As $G$ has no cycles, we place three closed loop poles at $-3$ and two at $-2$. The $\alpha$'s are unique upto the ordering and they are indeed the system poles. So we choose $\alpha_1=\alpha_3=\alpha_5=-3$ and $\alpha_2=\alpha_4=-2$. The formation is specified by the relative distance between agents so as to form a equilateral hexagon of side $2$. We first achieve a hexagon formation at leader position $f_n=(3,-1.829)$. Then the desired leader position is changed to $f_n'=(7,-3.829)$ instantly. We assume that the $\beta$'s can be computed in a time frame much faster than the smallest time constant determined by the $\al$'s and that the $\beta$'s are computed according to the algorithm given in \ref{subsec:rmk}. The resulting motion is demonstrated in figure \ref{fig:moving-form}
 \begin{figure}
  \begin{center}
    \includegraphics[width=0.3\textwidth]{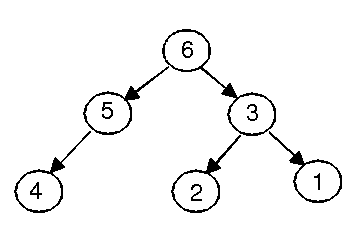}
  \end{center}
  \caption{Communication graph for example \ref{exmp:4}}\label{fig:graph-tree} 
\end{figure}
\begin{figure}
  \begin{center}
    \includegraphics[width=0.52\textwidth]{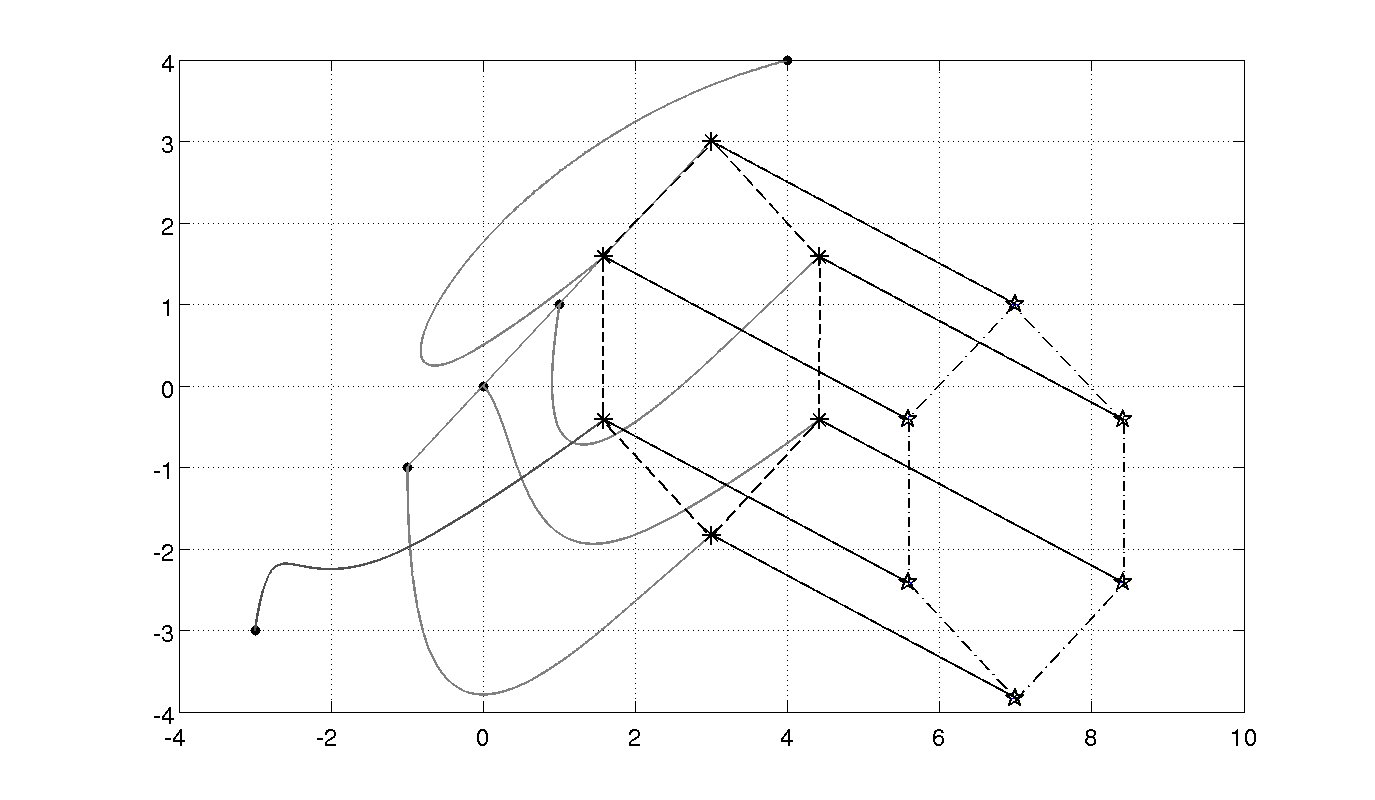}
  \end{center}
  \caption{Formation with changing $f_n$ by locally changing $\beta's$}\label{fig:moving-form} 
\end{figure} 
\end{example}

\section{Conclusion}
In this article we propose an algorithm for decentralized formation control which can achieve arbitrary pole placement for a multi-agent system in 
a leader-follower configuration. The  input
to leader is external and was required only to drive the leader to the desired state as 
specified in the formation. 
We have shown that if the directed communication graph is acyclic 
then the pole of each agent is determined by the self feedback gain. The same
condition proved to be sufficient for the ability to design the inter-agent gains that ensures rest of the agents follow the leader appropriately
to reach the formation. It would be interesting to prove that
directed cycles are helpful (and required) for
achieving arbitrary pole placement: both real and complex.

\bibliographystyle{ieeetr}
\bibliography{multiagent}

\begin{thebibliography}{10}

\bibitem{Ren2004}
W.~Ren and R.~Beard, {\em {Distributed consensus in mutivehicle cooperative
  control: Theory and Applicatons}}.
\newblock Springer, 2004.

\bibitem{Fax2004}
J.~Fax and R.~Murray, ``{Information Flow and Cooperative Control of Vehicle
  Formations},'' {\em IEEE Transactions on Automatic Control}, vol.~49,
  pp.~1465--1476, Sept. 2004.

\bibitem{Olfati2004}
R.~Olfati-Saber and R.~Murray, ``Consensus problems in networks of agents with
  switching topology and time-delays,'' {\em IEEE Transactions on Automatic
  Control}, vol.~49, no.~9, pp.~1520--1533, Sept.2004.

\bibitem{Olfati-Saber2007}
R.~Olfati-Saber, J.~Fax, and R.~Murray, ``{Consensus and cooperation in
  networked multi-agent systems},'' {\em Proceedings of the IEEE}, vol.~95,
  no.~1, pp.~215--233, 2007.

\bibitem{Jadbabaie2003}
A.~Jadbabaie and A.~Morse, ``{Coordination of groups of mobile autonomous
  agents using nearest neighbor rules},'' {\em IEEE Transactions on Automatic
  Control}, vol.~48, pp.~988--1001, June 2003.

\bibitem{Ren2005p}
W.~Ren and R.~Beard, ``Consensus seeking in multiagent systems under
  dynamically changing interaction topologies,'' {\em IEEE Transactions on
  Automatic Control}, vol.~50, no.~5, pp.~655--661, May,2005.

\bibitem{Fax2001}
J.~A. Fax, {\em Optimal and cooperative control of vehicle formations}.
\newblock PhD thesis, California Institute of Technology, 2001.

\bibitem{Lafferriere2005}
G.~Lafferriere, a.~Williams, J.~Caughman, and J.~Veerman, ``{Decentralized
  control of vehicle formations},'' {\em Systems \& Control Letters}, vol.~54,
  pp.~899--910, Sept. 2005.

\bibitem{lovasz}
L.~Lovasz and M.~Plummer, {\em Matching Theory}.
\newblock North Holland: Elsevier Science Publishers, 1986.

\bibitem{asratian}
A.~Asratian, T.~Denley, and R.~Haggkvist, {\em Bipartite Graphs and their
  Applications}.
\newblock United Kingdom: Cambridge University Press, 1998.

\bibitem{Chatterjee1977}
S.~Chatterjee and E.~Seneta, ``Towards consensus: Some convergence theorems on
  repeated averaging,'' {\em Journal of Applied Probability}, pp.~89--97, 1977.

\bibitem{Borkar1982}
V.~Borkar and P.~Varaiya, ``{Asymptotic agreement in distributed estimation},''
  {\em IEEE Transactions on Automatic Control}, pp.~0--5, 1982.

\bibitem{Hong2008}
Y.~Hong, G.~Chen, and L.~Bushnell, ``Technical communique: Distributed
  observers design for leader-following control of multi-agent networks,'' {\em
  Automatica (Journal of IFAC)}, vol.~44, no.~3, pp.~846--850, 2008.

\bibitem{Olfati-Saber2006}
R.~Olfati-Saber, ``Flocking for multi-agent dynamic systems: Algorithms and
  theory,'' {\em IEEE Transactions on Automatic Control}, vol.~51, no.~3,
  pp.~401--420, 2006.

\bibitem{Ren2008}
W.~Ren, ``On consensus algorithms for double-integrator dynamics,'' {\em IEEE
  Transactions on Automatic Control}, vol.~53, no.~6, pp.~1503--1509, 2008.

\bibitem{Ren2007}
W.~Ren, ``Second-order consensus algorithm with extensions to switching
  topologies and reference models,'' in {\em American Control Conference, 2007.
  ACC'07}, pp.~1431--1436, IEEE, 2007.

\bibitem{Yu2010}
W.~Yu, G.~Chen, and M.~Cao, ``Some necessary and sufficient conditions for
  second-order consensus in multi-agent dynamical systems,'' {\em Automatica},
  vol.~46, no.~6, pp.~1089--1095, 2010.

\bibitem{bab:fra}
L.~Babai and P.~Frankl, {\em Linear Algebraic Methods in Combinatorics}.
\newblock University of Chicago: Department of Computer Science, 1992.

\bibitem{RosenthalWillems}
J.~Rosenthal, M.~Schumacher, and J.~C. Willems, ``Generic eigenvalue assignment
  by memoryless real output feedback,'' {\em Systems and Control Letters},
  vol.~26, no.~4, pp.~253--260, 1995.

\end{thebibliography}
\end{document}